\documentclass[a4paper]{amsart}
\usepackage{amssymb}
\usepackage{stmaryrd}
\usepackage[all]{xy}
\usepackage{epsf}

\newcommand     {\printname}[1] {}
\newcommand{\labell}[1] {\label{#1}\printname{#1}}

\newtheorem{theorem}{Theorem}[section]
\newtheorem{proposition}[theorem]{Proposition}

\newtheorem{example}[theorem]{Example}
\newtheorem*{problem}{Linearization Problem}
\newtheorem*{problem1}{Open Problem}

\theoremstyle{remark}
\newtheorem{remark}[theorem]{Remark}

\numberwithin{equation}{section}

\newcommand{\Rr}{\mathbb R}

\newcommand{\Cc}{\mathbb C}

\newcommand{\Nn}{\mathbb N}

\newcommand{\set}[1]{\left\{#1\right\}}

\newcommand{\del}[2]{\frac{\partial #1}{\partial #2}}

\newcommand{\X}{\ensuremath{\mathfrak{X}}}
\newcommand{\M}{\ensuremath{\mathfrak{M}}}
\newcommand{\R}{\ensuremath{\mathfrak{R}}}

\newcommand{\tto}{\rightrightarrows}

\newcommand{\G}{\mathcal{G}}            
\newcommand{\s}{\mathbf{s}}             
\renewcommand{\t}{\mathbf{t}}           

\newcommand{\al}{\alpha}                
\newcommand{\be}{\beta}                 
\newcommand{\Lie}{\mathcal{L}}          
\renewcommand{\gg}{\mathfrak{g}}        
\newcommand{\hh}{\mathfrak{h}}        
\renewcommand{\ll}{\mathfrak{l}}        
\newcommand{\rr}{\mathfrak{r}}        
\renewcommand{\ss}{\mathfrak{s}}        
\DeclareMathOperator{\Ker}{Ker}         
\DeclareMathOperator{\End}{End}         

\newcommand{\comment}[1]{}

\begin{document}

\title[Linearization]{Linearization of Poisson brackets}

\author{Rui Loja Fernandes}
\address{Depart.~de Matem\'{a}tica, Instituto Superior T\'{e}cnico,
  1049-001 Lisboa, PORTUGAL}
\email{rfern@math.ist.utl.pt}
\thanks{Supported in part by FCT/POCTI/FEDER and by grant
  POCTI/1999/MAT/33081.}

\author{Philippe Monnier}
\address{Center for Mathematical Analysis, Geometry and Dynamical
  Systems, Instituto Superior T\'{e}cnico,
  1049-001 Lisboa, PORTUGAL}
\email{pmonnier@math.ist.utl.pt}

\subjclass[2000]{}

\date{February, 2004}

\begin{abstract}
We review the linearization of Poisson brackets and related problems,
in the formal, analytic and smooth categories.
\end{abstract}
\maketitle

\section{Introduction}             %
\labell{Linearization: Outline}    %

Recall that a \textbf{Poisson bracket} on a manifold $M$ is a Lie bracket
$\{\cdot,\cdot\}:C^{\infty}(M)\times C^{\infty}(M)\to C^{\infty}(M)$,
satisfying the derivation property
\[ \{ fg, h\}= f\{ g, h\}+ g\{ f, h\}, \quad f,g,h \in C^{\infty}(M).\]
The Weinstein Splitting Theorem (see \cite{Wein}, Theorem 2.1) states
that a neighborhood of a point $x_0\in M$ is Poisson diffeomorphic to
the direct product of a symplectic manifold and a Poisson manifold for
which the bracket vanishes at $x_0$. Since a symplectic manifold has
no local invariants (Darboux's Theorem), the local study of Poisson
brackets reduces to the study of brackets vanishing at $x_0$.

Let $\{\cdot,\cdot\}$ be a Poisson bracket in $M$ that vanishes at $x_0$.
If we choose local coordinates $(x^1,\dots,x^m)$ in a neighborhood $U$ 
of $x_0$, the Taylor series of the function $\{x^i,x^j\}$ reads
\[ \{x^i,x^j\}(x)=c^{ij}_k x^k+g^{ij}_k(x)x^k,\]
where
$c^{ij}_k=\del{\{x^i,x^j\}}{x^k}(x_0)$,
and the $g^{ij}_k$ are smooth functions vanishing at $x_0$. In this
paper we are interested in the following \emph{linearization problem}:

\begin{problem}
  Are there new coordinates where the functions $g^{ij}_k$ vanish
  identically, so that the bracket is linear in these coordinates?
\end{problem}

More generally, we will be interested in local canonical forms for
the Poisson bracket.

Alternative algebraic and geometric descriptions of this problem will
be given below. The constants $c^{ij}_k$ are the structure constants
of a Lie algebra $\gg$, called the \textbf{isotropy Lie algebra} at
$x_0$. It can be defined independently of the choice of coordinates. 
Both the algebraic and the geometric descriptions will show the 
relevance of this Lie algebra to the linearization problem. 

The linearization problem can be studied in different setups. We
have just described it in the smooth category, which is the most
challenging version of the problem. We are also interested in the
following categories:
\begin{itemize}
\item \emph{Formal category}: we consider a Poisson bracket on a
  formal neighborhood of the origin, and we look for a formal change
  of coordinates where the functions $g^{ij}_k$ vanish to all orders at $x_0$.
\item \emph{Analytic category}: we consider an analytic Poisson bracket
  on an analytic manifold, and we look for an analytic change of coordinates
  where the functions $g^{ij}_k$ vanish.
\end{itemize}

The first results on linearization of Poisson brackets are due to
V.I. Arnold (\cite{Arnold}) and Alan Weinstein. In the foundational 
paper \cite{Wein}, Weinstein showed that
the formal linearization problem can be reduced to a cohomology
question which can always be solved if the isotropy Lie algebra is
\emph{semisimple}. For the other categories, he did not give definite
results, but he made the following conjectures:
\begin{enumerate}
\item[(i)] for a semisimple isotropy Lie algebra, analytic
  linearization of Poisson brackets can always be achieved.
\item[(ii)] for a semisimple isotropy Lie algebra of compact type, 
  smooth linearization of Poisson brackets can always be achieved.
\end{enumerate} 
Weinstein also gave a counter-example of a smooth Poisson bracket with
isotropy Lie algebra $\mathfrak{sl}(2)$ which is not linearizable.

Both Weinstein conjectures were proved later by Jack Conn in
\cite{Conn1} and \cite{Conn2}. Conn's proofs are hard analysis
proofs. They consist in showing that the formal power series change of
coordinates obtained by Weinstein, which formally linearizes the
Poisson bracket, actually yield linearization coordinates. For the
analytic case, Conn uses methods due to Arnol'd and Kolmogoroff to
prove convergence, while for the smooth case, he applies the
Nash-Moser approximation scheme, used by Nash on his famous proof of
the isometric embedding theorem, to obtain convergence. We will outline
the main steps of the proofs below.

After Conn's work was completed, attention turned to other type of
isotropy Lie algebras. Let us call a Lie algebra $\gg$
\textbf{non-degenerate relative to Poisson brackets} if any Poisson
structure with isotropy Lie algebra isomorphic to $\gg$ is
linearizable, and \textbf{degenerate} otherwise. In
\cite{Wein:rrank2}, Weinstein showed that semisimple Lie algebras of
real rank greater than one are degenerate. The case of real rank 1,
with the exception of $\mathfrak{sl}(2,\Rr)$, remains open. In
\cite{Dufour:Poisson}, Dufour has classified all non-degenerate
3-dimensional Lie algebras. More recently, in \cite{JP-Zung} Dufour
and Zung have proved formal and analytic non-degeneracy for the Lie
algebra of affine transformations $\mathfrak{aff}(n)$. On the other
hand, in \cite{BaCr}, examples are given of Poisson structures for
which linearization can be decided only from knowledge of its higher
order jets. The ultimate goal, which seems beyond our current state of
knowledge, would be:
\begin{problem1}
  Characterize the non-degenerate Lie algebras relative to Poisson brackets.
\end{problem1}
For linearization around singular symplectic leaves, we refer the reader 
to the recent preprint \cite{Bra}.

In this paper we will survey what is known about the linearization of
Poisson brackets. We will focus on comparing the analytic approach,
which one can find in the literature, and the geometric approach which
is the subject of recent work (see \cite{CrFe2,CrFeMo,Zung2}), and
which sheds much light on this problem. During the preparation of this
survey we were informed by Dufour and Zung of their upcoming book
\cite{Book:JP-Zung} on normal forms for Poisson structures, where the
material discussed here (and much more) will be presented in greater
detail.

The problem of linearization of Poisson brackets has remarkable
similarities with the problem of linearization of Lie algebra
actions. In Section \ref{sec:Lie:algebras} below, we will recall the
classical results of Cartan (\cite{Cartan}), Bochner (\cite{Bochner}),
Hermann (\cite{Hermann}) and Guillemin and Sternberg (\cite{GuiSter}),
on linearization of Lie group and Lie algebra actions, in the formal,
analytic and smooth categories. This easier problem sets up the stage
for the study of the linearization of Poisson brackets. Again, we will
see here the advantages of the geometric approach over the analytic
approach.

We take up the linearization of Poisson brackets in Section
\ref{sec:Poisson:brackets} below. We discuss the formal, analytic and
smooth linearization problems, and we sketch the known proofs (without
any technical details). These are all analytic proofs. We will mention
briefly the geometric proofs proposed in \cite{CrFe2,CrFeMo,Zung2}
(analytic case) and in \cite{CrFe2} (smooth case). We shall also
emphasize the similarity with the linearization of Lie algebra
actions, though the known proofs are quite different in each case. 

In section 4, we turn to a more general normal form: the so called
\emph{Levi decomposition} for a Poisson bracket. This is an analogue
for Poisson brackets of the usual Levi decomposition for finite
dimensional Lie algebras. After all, Poisson brackets are infinite
dimensional Lie brackets, and so one should expect some kind of Levi
decomposition to hold as well. This Levi decomposition may also be seen
as a semi-linearization of Poisson brackets. Again, there is a formal
Levi decomposition due to Wade (\cite{Wade}), an analytic Levi
decomposition due to Zung (\cite{Zung}), and a smooth Levi decomposition
due to Monnier and Zung (\cite{MoZung}).

Finally, in Section \ref{sec:algebroids}, we give a short overview of
linearization and Levi decomposition in the context of Lie algebroids
which allows one to unify the results for actions of Lie algebras and
for Poisson brackets.

We would like to thank Marius Crainic and Nguyen Tien Zung for several 
discussions related to linearization of Poisson brackets, and Alan Weinstein 
for comments on an earlier version of this manuscript.
 
\section{Linearization of Lie algebra actions}
\labell{sec:Lie:algebras}

The problem of linearization of Poisson brackets has remarkable
similarities with the problem of linearization of Lie algebra
actions. First of all, for group actions, we have the following
well-known classical result:

\begin{theorem}[Bochner \cite{Bochner}]
Let $G\times M\to M$ be a smooth action of a compact Lie group on a
smooth manifold, which has a fixed point. Then, in suitably chosen local
coordinates around the fixed point, the action is linear.
\end{theorem}

For a modern version of the proof, we refer to the book by Duistermaat
and Kolk (\cite{DuKo}). This result grew out of H.~Cartan's investigations
on analytic group actions, and Cartan himself proved a similar result
in the analytic case (\cite{Cartan}). Palais and Smale suggested
extending this result to non-compact Lie group actions and Hermann
(\cite{Hermann}) explained that the corresponding formal problem can
be reduced to a cohomology question which can always be solved if the
group is semisimple. 

If the group $G$ is connected, trying to find a linear system of
coordinates for the Lie group action, is the same as finding a linear
system of coordinates for the associated Lie algebra action. Note,
however, that since a Lie algebra action does not always integrate to
a Lie group action, the linearization problem for Lie algebra actions
is harder. Let us look then at this problem.

Consider an action $\rho:\gg\to\X(M)$ of a Lie algebra $\gg$ on a
manifold $M$, which has a fixed point $x_0\in M$. We are interested in
the following problem.

\begin{itemize}
\item Are there coordinates $(x^1,\dots,x^m)$ around $x_0$ such that,
  for all $X\in\gg$, the vector fields $\rho(X)$ are linear in these
  coordinates ?
\end{itemize}

In general, if we fix a basis $\{X_1,\dots,X_N\}$ for $\gg$ and an
arbitrary system of coordinates $(x^1,\dots,x^m)$ around $x_0$, the
action will be written in the form\footnote{We will make systematic
use of the Einstein sum convention.}:
\begin{equation}
\labell{eq:action}
 \rho(X_\al)=c_{\al j}^i x^j \del{}{x^i}+O(2), \qquad (\al=1,\dots, N),
\end{equation}
and so we are looking for coordinates where the higher order terms
vanish.

In order to put this problem in a more geometric framework, observe that the
action $\rho$ induces a linear action $\rho_L$ on the tangent space
$T_{x_0} M$. In terms of a basis for $\gg$ and a system of coordinates
around $x_0$ as above, the linear action is just given by
truncating in (\ref{eq:action}) the higher order terms: If
$(u^1,\dots,u^m)$ denotes the linear coordinates defined by the basis
$\{\left.\del{}{x^1}\right|_{x_0},\dots,\left.\del{}{x^m}\right|_{x_0}\}$
for $T_{x_0} M$, the linear action is defined by
\begin{equation}
\labell{eq:linear:action}
 \rho_L(X_\al)=c_{\al j}^i u^j \del{}{u^i}\qquad (\al=1,\dots, N).
\end{equation}

The \emph{geometric formulation} of the linearization problem is then:
\begin{itemize}
\item Is there a local diffeomorphism $\phi:M\to T_{x_0}M$, from a
  neighborhood of $x_0$ in $M$ to a neighborhood of $0$ in
  $T_{x_0}M$, which is $\gg$-equivariant?
\end{itemize}

Let us denote by $\M$ the maximal ideal formed by germs of smooth functions
that vanish at $x_0$:
\[ \M=\set{ f\in C^{\infty}(M): f(x_0)=0 }.\]
There is also an induced linear action, denoted by $\rho_L^*$, of $\gg$ on
$\M$:
\[ \rho_L^*(X)(f)=X(f).\]
It is easy to see that $\M^2\subset \M$ is an invariant ideal, so both
$\M^2$ and the cotangent space $T^*_{x_0}M\simeq \M/\M^2$ have induced
linear $\gg$-actions. Obviously, the $\gg$-actions on $T^*_{x_0}M$ and
$T_{x_0}M$ are transpose to each other. The \emph{algebraic formulation}
of the linearization problem is:

\begin{itemize}
\item Is there a splitting $\sigma$ of the short exact sequence of
  $\gg$-modules
  \begin{equation}
    \labell{eq:seq:Lie:module}
    \xymatrix{ 0\ar[r]& \M^2\ar[r]& \M\ar[r]&
      \M/\M^2\ar[r] \ar@(dl,dr)@{-->}[l]^{\sigma}&0}\ ?
  \end{equation}
\end{itemize}
In fact, the projection $\M\to \M/\M^2$ is just taking the
differential at $x_0$: $f\mapsto d_{x_0}f$. Hence, if
$\sigma:\M/\M^2\to\M$ is such a splitting and we fix a basis
$\{\xi^1,\dots,\xi^n\}$ for $\M/\M^2$, then $\sigma(\xi^i)$ are germs
of smooth functions $x^1,\dots,x^m$, for which the differentials
$d_{x_0}x^i=\xi^i$ are independent. Therefore, $(x^1,\dots,x^m)$
yields a system of coordinates around $x_0$. On the other hand, since
$\sigma$ is $\gg$-equivariant, when we express the action in this
system of coordinates we obtain a linear action.

The linearization problem we have just explained was set up in the
smooth category. Similarly, we can consider formal linearization
(replace $C^\infty(M)$ by the ring of power series
$\Rr[[x^1,\dots,x^m]]$) or analytic linearization (replace
$C^\infty(M)$ by the ring of analytic functions $C^\omega(M)$).

\subsection{Formal Linearization}
The algebraic version of the linearization problem suggests one should
look at obstructions coming from the Lie algebra cohomology of
$\gg$. These are the only obstructions for the formal problem, as we
show now.

Recall that if $V$ is any $\gg$-module, then we have the
Chevalley-Eilenberg complex $(C^\bullet(\gg;V),d)$, where: 
\begin{itemize}
\item $C^r(\gg;V)$ is the vector space consisting of all
  $r$-multilinear alternating forms
  $\omega:\gg\wedge\cdots\wedge\gg\to V$, and 
\item the differential $d:C^r(\gg;V)\to C^{r+1}(\gg;V)$ is defined by
\begin{multline}
\labell{eq:differential}
d\omega(X_0,\dots,X_r)=\sum_{i=0}^{r}
(-1)^k X_k(\omega(X_0,\dots,\widehat{X}_i,\dots,X_r))\\
+\sum_{i<j}(-1)^{i+j}
\omega([X_i,X_j],X_0,\dots,\widehat{X}_i,\dots,\widehat{X}_j,\dots,X_r).
\end{multline}
\end{itemize}
The cohomology of this complex is called the \textbf{Lie algebra cohomology}
of the $\gg$-module $V$, and will be denoted $H^\bullet(\gg;V)$.

It is a well-known fact from Homological Algebra (see, e.g., \cite{Weibel})
that the obstructions to the existence of a splitting of a short exact
sequence of $\gg$-modules lie in the first Lie algebra cohomology. Let
us recall this briefly. Suppose we are given a short exact sequence of
$\gg$-modules
\[ \xymatrix{0\ar[r]& V\ar[r]& W\ar[r]^{\pi}&Z\ar[r]&0}.\]
Let $\sigma:Z\to W$ be a vector space splitting for this sequence. The
space of linear endomorphisms $\End(Z,V)$ is naturally a $\gg$-module, and
we can define an element $\omega\in C^1(\gg;\End(Z,V))$ by:
\[ \omega(X): z\longmapsto \sigma(X\cdot z)-X\cdot\sigma(z)\qquad
(z\in Z).\]
Notice that the left-hand side actually lies in $\Ker\pi=V$. A
quick computation shows that $d\omega=0$, so we have a well-defined
class $[\omega]\in H^1(\gg;\End(Z,V))$. It is easy to check that
(i) this class is independent of the splitting, and (ii) it vanishes iff
there exists a splitting of the sequence, as $\gg$-modules.

Now, the First Whitehead Lemma states that if $\gg$ is a semisimple
Lie algebra and $V$ is any \emph{finite dimensional} $\gg$-module, then
$H^1(\gg;V)=0$. Certainly, we cannot apply this directly to the short
exact sequence (\ref{eq:seq:Lie:module}), since the $\gg$-modules are
infinite dimensional. Still, one has:

\begin{theorem}[Hermann \cite{Hermann}]
\labell{thm:Hermann}
If $\gg$ is a semisimple Lie algebra, then the action is formally
linearizable, i.e., there exist coordinates
$(x^1_\infty,\dots,x^m_\infty)$ around $x_0$ such that, for any basis
$\{X_1,\dots,X_N\}$ for $\gg$, we have
\[  \rho(X_\al)=c_{\al j}^i x^j_\infty \del{}{x^i_\infty}+o(\infty),
\qquad (\al=1,\dots, N),\]
where $c_{\al j}^i$ are some constants.
\end{theorem}

Let us sketch a proof of this result. Instead of considering the full
ring $\M$ at once, one looks at the spaces $\M/\M^{k+1}$ (germs mod
terms of order higher than $k$). Notice that we still have short exact
sequences of $\gg$-modules
\begin{equation}
\labell{eq:seq:aux:Lie:algebra}
\xymatrix{
0\ar[r]& \M^k/\M^{k+1}\ar[r]&\M/\M^{k+1}\ar[r]^{\pi_k}&\M/\M^k\ar[r]&0
}.
\end{equation}
In this sequence, all $\gg$-modules are finite dimensional. Hence,
by the First Whitehead Lemma, we have splittings
$\phi_k:\M/\M^k\to \M/\M^{k+1}$. It follows that one has a
commutative diagram
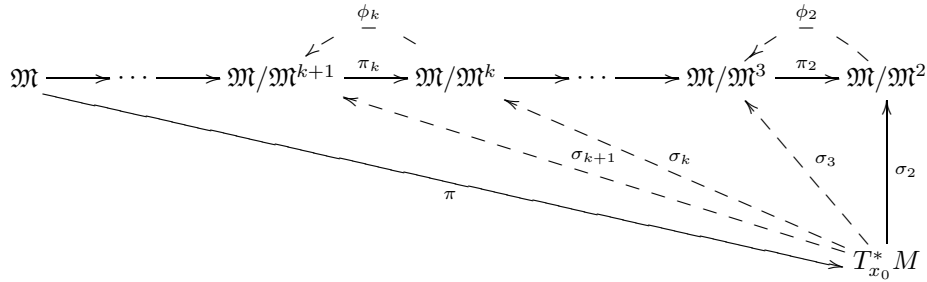
\[
\xymatrix{
\M\ar[r]\ar@<-1ex>[ddrrrrrr]_{\pi}&
\cdots\ar[r]&\M/\M^{k+1}\ar[r]^{\pi_k}
&\M/\M^k\ar[r]\ar@(ul,ur)@{-->}[l]_{\phi_k}&\cdots\ar[r]&
\M/\M^3\ar[r]^{\pi_2}&\M/\M^2\ar@(ul,ur)@{-->}[l]_{\phi_2}\\
\\
&&&&&&T^*_{x_0}M\ar[uu]_{\sigma_2}
\ar@{-->}[uul]_{\sigma_3}\ar@{-->}[uulll]_{\sigma_k}
\ar@{-->}[uullll]_{\sigma_{k+1}}}
\]
In this diagram, $\sigma_2:T^*_{x_0}M\to \M/\M^2$ denotes the natural
isomorphism, while for $k>2$ the $\gg$-module homomorphisms
$\sigma_k:T^*_{x_0}M\to\M/\M^k$ are obtained by induction:
\[ \sigma_k=\phi_k\circ\sigma_{k-1}.\]
By passing to the projective limit, we obtain a $\gg$-module homomorphism
\[ \sigma: T^*_{x_0} M\to \varprojlim \M/\M^k\subset \M,\]
and, as above, this defines the coordinate system
$(x^1_\infty,\dots,x^m_\infty)$ that formally linearizes the action.

\subsection{Analytic Linearization}

The proof given above for Hermann's Linearization Theorem shows that
one constructs the formal linearizing coordinates by a succession of changes
of coordinates, which at each step removes terms of order $k$, without
modifying the terms of order less than $k$. Hence, it is natural to
ask if one can make the choices at each step so that, in the limit, one
obtains a smooth or analytic coordinate system, that linearizes the
action. In the analytic category this is indeed the case:

\begin{theorem}[Guillemin and Sternberg \cite{GuiSter}]
\labell{thm:lie:analytic}
If $\gg$ is a semisimple Lie algebra
and the action is analytic, then it is analytically linearizable,
i.e., there exist analytic coordinates
$(x^1_\infty,\dots,x^m_\infty)$ around $x_0$ such that, for any
basis $\{X_1,\dots,X_N\}$ for $\gg$, we have
\[  \rho(X_\al)=c_{\al j}^i x^j_\infty \del{}{x^i_\infty},
\qquad (\al=1,\dots, N),\] for some constants $c_{\al j}^i$.
\end{theorem}

\begin{proof}[Analytic proof]
One proves, directly, convergence of the formal linearizing
coordinates. The estimates are rather involved, and so we will omit
them. The interested reader can obtain them by specializing the proof
of the analytic Levi decomposition for Lie algebroids, given recently
in \cite{Zung}, to the case of the action Lie algebroid.
\end{proof}

\begin{proof}[Geometric proof]
This is the proof given by Guillemin and Sternberg in \cite{GuiSter}. 
It avoids the question of convergence of the formal linearizing
coordinates, through the use of complexification and analytic
continuation into the complex plane, to obtain an action of the
complexified Lie algebra. Restricting this action to the compact real
form, one obtains an action which integrates to an action of a
compact Lie group. This was proved by Guillemin and Sternberg using
the assumption that the action is analytic, but in the next section we
will show that this is already true for smooth actions. We can then
apply Bochner's Theorem to linearize the action.
\end{proof}

\begin{remark}
Independently of Guillemin and Sternberg, but using similar
techniques, Kushnirenko (\cite{Kushnirenko}) proved linearization of
real-analytic or complex-analytic actions of local real or complex
semisimple Lie groups.
\end{remark}

\begin{remark}
\label{rem:flato}
A completely different proof was obtained a few years after
\cite{GuiSter}, by Flato, Pinczon and Simon (\cite{Flato}), as a
corollary of their study of nonlinear actions of Lie groups and Lie
algebras having a fixed point. Their main objective was the infinite
dimensional case needed to study the linearizability of nonlinear
evolution equations covariant under a Lie group action. However, the
Lie theory that was developed (the Lie group-Lie algebra analytic
nonlinear representations correspondence), combined with Weyl's
unitary trick to bring the problem to compact forms, and a variant of
Bochner's theorem above, leads to a simple proof of \ref{thm:lie:analytic}.
\end{remark}

\subsection{Smooth Linearization}

The following example of a smooth action of $\mathfrak{sl}(2,\Rr)$
which is not linearizable is due to Guillemin and Sternberg \cite{GuiSter}.

\begin{example}
\labell{ex:non:linear:action}
Consider the basis $\{X,Y,Z\}$ of $\mathfrak{sl}(2,\Rr)$ satisfying
the relations:
\[ [X,Y]=-Z,\quad [Y,Z]=X, \quad [Z,X]=Y.\]
We have a linear action defined by:
\begin{align*} 
\rho_L(X)=y\frac{\partial}{\partial z}+z\frac{\partial}{\partial y},&\quad
\rho_L(Y)=x\frac{\partial}{\partial z}+z\frac{\partial}{\partial x},\\
\rho_L(Z)=&x\frac{\partial}{\partial y}-y\frac{\partial}{\partial x}.
\end{align*}
For this action the orbits are the level sets of the quadratic form
$x^2+y^2-z^2=r^2-z^2$. On the other hand, we can perturb this action to the
non-linear action:
\begin{align*} 
\rho(X)=\rho_L(X)+\frac{xz}{r^2}g(r^2-z^2)V,&\quad
\rho(Y)=\rho_L(Y)-\frac{yz}{r^2}g(r^2-z^2)V,\\
\rho(Z)=\rho_L(Z)&+g(r^2-z^2)V,
\end{align*}
where $V=x\frac{\partial}{\partial x}+y\frac{\partial}{\partial y}+
z\frac{\partial}{\partial z}$ is the radial vector field,
and $g\in C^\infty(\Rr)$ is such that $g(x)>0$, if $x>0$, and $g(x)=0$,
if $x\le 0$.

The orbits of $\rho$ coincide with the orbits of $\rho$ inside the
cone $r=z$. Outside this cone, the orbits of $\rho(Z)$ spiral towards
the origin while the orbits of $\rho_L(Z)$ are circles. Hence, $\rho$
is not linearizable.
\end{example}

This shows that, in the smooth case, we need to restrict the class of
semisimple Lie algebras. Recall that a finite dimensional Lie algebra $\gg$ is 
called \textbf{semisimple of compact type} if any of the following
equivalent conditions are satisfied (see, e.g., 
\cite{DuKo}, Section 3.6):
\begin{enumerate}
\item[(i)] The Killing form of $\gg$ is negative definite;
\item[(ii)] The simply connected Lie group integrating $\gg$ is
  compact;
\item[(iii)] Any Lie group integrating $\gg$ is compact;
\end{enumerate}
For such Lie algebras we have the following linearization result:

\begin{theorem}
\labell{thm:infin:Bochner}
If $\gg$ is a semisimple Lie algebra of compact type and the action is
smooth, then it is smoothly linearizable, i.e., there exist smooth
coordinates $(x^1_\infty,\dots,x^m_\infty)$ around $x_0$ such that,
for any basis $\{X_1,\dots,X_N\}$ for $\gg$, we have
\[  \rho(X_\al)=c_{\al j}^i x^j_\infty \del{}{x^i_\infty},
\qquad (\al=1,\dots, N),\]
for some constants $c_{\al j}^i$.
\end{theorem}

To our surprise, in spite of Guillemin and Sternberg work being more
that 30 years old, we could not find this result in the literature. 

\begin{proof}[Analytic proof]
Again, one proves, directly, convergence of the formal linearizing
coordinates, and the estimates are even more involved than in the
analytic case. Once more, they can also be obtained by specializing the proof
of the smooth Levi decomposition for Lie algebroids, given recently
in \cite{MoZung}, to the case of the action Lie algebroid (see Section
\ref{sec:Levi}).
\end{proof}

\begin{proof}[Geometric proof]
The proof consists in showing that the $\gg$-action can be integrated
to a $G$-action, where $G$ is the compact simply connected Lie group
integrating $\gg$, so we can apply Bochner's Linearization Theorem
for compact Lie group actions.

A Lie algebra action $\rho:\gg\to\X(M)$ in general will not integrate
to a Lie group action. Our approach is to search for a Lie group
action by looking at its graph. This graph should be obtained by integrating
the graph of the Lie algebra action. Since the graph of a Lie group
(respectively, Lie algebra) action is in fact the associated action
Lie groupoid (resp.~Lie algebroid) what is involved here is the
integration of Lie algebroids to Lie groupoids.

For \emph{any} Lie algebra action $\rho:\gg\to\X(M)$ the associated
action Lie algebroid always integrate to a Lie groupoid. This fact was
first proved by Dazord (\cite{MoMr} and \cite{CaWe}, Section 16.4) and
follows also easily from the general integrability result in
\cite{CrFe1}. Let $\s,\t:\G\rightrightarrows M$ be the Lie groupoid
integrating the action Lie algebroid.
Since $x_0$ is a fixed point for the action, the source and target
satisfy $\s^{-1}(x_0)=\t^{-1}(x_0)=G$, a simply connected Lie group
integrating $\gg$. So this fiber is compact, and by stability, we
conclude that $\G$ trivializes over an neighborhood $U$ of $x_0$ with
compact closure. Since this neighborhood is invariant for the action,
this shows that the action in $U$ integrates. 
\end{proof}

\begin{remark}
In general, a Lie groupoid is not assumed to be Hausdorff, but its $s$-fibers
are. It can be shown that, in our case, the Lie groupoid that integrates the
action is Hausdorff, so one can apply the usual Reeb stability
result. Alternatively, one can use the stability result for
non-Hausdorff manifolds due to Mr\v{c}un (see \cite{Mr}). 
\end{remark}
\vskip 10 pt

Of course, just like in the Poisson case, one can ask for which Lie
algebras $\gg$ is it true that every $\gg$-action is linearizable
around a fixed point. Call such a Lie algebra non-degenerate relative
to actions.

\begin{problem1}
Characterize the Lie algebras that are non-degenerate relative to
actions. 
\end{problem1}

The $\mathfrak{sl}(2,\Rr)$-action on $\Rr^3$ given by Gullemin and Sternberg,
does not integrate to an $SL(2,\Rr)$-action. In \cite{CairGhy}, Cairns and 
Ghys modify this example to obtain an $SL(2,\Rr)$-action on $\Rr^3$ that 
is not linearizable. On the other hand, they show that a smooth 
$SL(n,\Rr)$-action on $\Rr^n$ is always linearizable. This suggests one 
should look for $\gg$ non-degeneracy for a fixed dimension of the manifold, or
even specifying the linear part of the action.

It is an intriguing question how this problem is related to the
question of non-degeneracy of the Lie algebra relative to Poisson
brackets. 

\section{Linearization of Poisson brackets}
\labell{sec:Poisson:brackets}

Let us now turn to Poisson geometry. Here we view a Poisson manifold
as an infinitesimal object. The corresponding global object (if it
exists) is a \emph{symplectic groupoid} $\Sigma\tto M$, and we can ask
for the analogue of Bochner's theorem. There is, in fact, such a
theorem: 

\begin{theorem}
\labell{thm:proper:groupoid}
Let $\Sigma\tto M$ be a symplectic groupoid with a fixed point $x_0\in
M$. Assume that the isotropy group $G_{x_0}=s^{-1}(x_0)=t^{-1}(x_0)$
is compact. Then $\Sigma$ is locally isomorphic to the symplectic
groupoid $T^*G_{x_0}\tto T_{x_0}M$. 
\end{theorem}

This result follows from a more general result valid for proper Lie
groupoids (see Weinstein \cite{Wein:proper}, Zung \cite{Zung3}). The
proof given in Zung \cite{Zung3} is analytic, and uses a fixed point
theorem. However, a much more geometric proof of this result, using
averaging (in fact, vanishing of cohomology \cite{Cr}) can be found in
\cite{CrFe2,CrFeMo}. Note that the isomorphism given by the theorem
preserves both the Lie groupoid structure and the symplectic
structure. If one ignores the Lie groupoid structure, then the
symplectic isomorphism follows from the usual Lagrangian Neighborhood
theorem of Weinstein, since the fiber $s^{-1}(x_0)=t^{-1}(x_0)$ over a
fixed point is always a Lagrangian submanifold.

Note, also, that the local isomorphism $\Sigma\simeq T^*G$ covers
a local Poisson diffeomorphism $M\simeq T_{x_0}M=\gg^*$ which maps a
neighborhood of $x_0\in M$ onto a neighborhood of $0\in\gg^*$. In this
way, trying to find a local isomorphism $\Sigma\simeq T^*G$ is the
same as trying to linearize the Poisson bracket. However, Poisson
manifolds do not always integrate to symplectic groupoids, and so the
problem of linearizing Poisson brackets is harder.  Let us look then
at this problem.

We consider a Poisson manifold $(M,\{~,~\})$ and we assume that the
bracket vanishes at $x_0\in M$.  If we fix local coordinates
$(x^1,\dots,x^m)$ around $x_0$, we have
\begin{equation}
\labell{eq:Poisson:brackets}
\{x^i,x^j\}=c^{ij}_k x^k+O(2).
\end{equation}
Here the $c^{ij}_k$ are structure constants for a Lie algebra, thus
the linear terms define a linear Poisson bracket. The linearization
problem we are interested in, is then:
\begin{itemize}
\item Are there coordinates $(x^1,\dots,x^m)$ around $x_0$ such that
  the Poisson bracket is linear in these coordinates?
\end{itemize}

In order to put this problem in a more geometric framework, recall
that the Poisson bracket gives rise to a Lie bracket $[\cdot, \cdot]$
on the space $\Omega^1(M)$ of 1-forms on $M$, which is given by:
\begin{equation}
  \labell{kozul:bracket}
  [\al, \be]=\Lie_{\#\al}\be-\Lie_{\#\be}\al-d\Pi(\al,\be),\quad
  \al,\be\in\Omega^1(M).
\end{equation}
Here $\Pi\in \Gamma(\Lambda^2TM)$ denotes the Poisson 2-tensor which
is associated to the Poisson bracket by $\Pi(df, dg)= \{f, g\}$, and
$\#:T^*M\to TM$ denotes contraction by $\Pi$.  This bracket satisfies
(in fact, it is uniquely determined by) the following two basic
properties:
\begin{enumerate}
\item[(i)] for exact 1-forms it coincides with the Poisson bracket:
  \[ [df, dg]= d\{f, g\},\quad f,g\in C^{\infty}(M);\]
\item[(ii)] it satisfies the Leibniz identity:
  \[ [\al, f\be]= f[\al, \be] + \#\be(f)\al, \quad
    \al,\be\in\Omega^1(M),\ f\in C^{\infty}(M).\]
\end{enumerate}
The triple $(T^*M,[\cdot, \cdot],\#)$ is a Lie algebroid, called
the \textbf{cotangent Lie algebroid} of the Poisson manifold $M$.

For each $x\in M$, the kernel $\gg_x=\Ker\#_x\subset T^*_x M$ is a Lie
algebra, called the \textbf{isotropy Lie algebra} at $x$. At $x_0$, we
have $\gg_{x_0}=T^*_{x_0}M$, so the dual $T_{x_0}M$ carries a Poisson
structure called the \textbf{linear approximation} at $x_0$. In local
coordinates $(x^1,\dots,x^m)$ around $x_0$, the numbers $c^{ij}_k$ in
(\ref{eq:Poisson:brackets}) are just the structure constants for
$\gg_{x_0}$ relative to the basis
$\{\left.dx^1\right|_{x_0},\dots,\left.dx^m\right|_{x_0}\}$.
The geometric reformulation of the linearization problem is then:
\begin{itemize}
\item Is there a local Poisson diffeomorphism $\phi:M\to T_{x_0}M$, from a
  neighborhood of $x_0$ in $M$ to a neighborhood of $0$ in
  $T_{x_0}M$?
\end{itemize}

Let us denote again by $\M$ the germs of smooth functions that vanish
at $x_0$. Then $\M$ is a maximal ideal of the Lie algebra
$(C^\infty(M),\{~,~\})$\footnote{Notice that, for Lie algebra actions,
  $\M$ was a \emph{$\gg$-module}, while now it is a \emph{Lie algebra}.}. 
The quotient $\M/\M^2=T^*_{x_0} M$ is just the isotropy Lie algebra
$\gg_{x_0}$, so we can reformulate the linearization problem in the following
algebraic form: 
\begin{itemize}
\item Is there a splitting $\sigma$ of the short exact sequence of
  Lie algebras
  \begin{equation}
    \labell{eq:seq:Lie:algebra}
    \xymatrix{ 0\ar[r]& \M^2\ar[r]& \M\ar[r]&
      \gg_{x_0}\ar[r] \ar@(dl,dr)@{-->}[l]^{\sigma}&0}\ ?
  \end{equation}
\end{itemize}
Again, the projection $\M\to \gg_{x_0}$ is just taking the
differential: $f\mapsto d_{x_0}f$. Hence, if $\sigma:\gg_{x_0}\to\M$
is such a splitting and we fix a basis $\{\xi^1,\dots,\xi^n\}$ for
$\gg_{x_0}$, then $\sigma(\xi^i)$ are germs of smooth functions
$x^1,\dots,x^m$, for which the differentials $d_{x_0}x^i=\xi^i$ are
independent. Therefore, $(x^1,\dots,x^m)$ yields a system of
coordinates around $x_0$. On the other hand, since $\sigma$ is a
homomorphism of Lie algebras, when we express the Poisson bracket in
this coordinate system we obtain a linear Poisson bracket.

Like the Lie algebra case, besides the smooth category, the
linearization problem for Poisson brackets can also be considered in
the formal category (replace $C^\infty(M)$ by the ring of power series
$\Rr[[x^1,\dots,x^m]]$) or the analytic category (replace
$C^\infty(M)$ by the ring of analytic functions $C^\omega(M)$).

\subsection{Formal Linearization}

Just like the case of Lie algebra actions, the algebraic version of
the linearization problem for Poisson brackets leads to obstructions
to linearization lying in Lie algebra cohomology.

In this case, we shall need to look at Lie algebra cohomology in
degree 2. In fact, it is a well-known fact in Homological Algebra (see,
e.g., \cite{Weibel}) that the obstructions to the existence of central
extensions lies in $H^2$. More exactly, let us consider a short exact
sequence of Lie algebras
\[ \xymatrix{0\ar[r]& V\ar[r]& \hh\ar[r]^{\pi}&\gg\ar[r]&0},\]
where $V$ is abelian (a vector space). If we fix a splitting
$\sigma:\gg\to \hh$ of this sequence as vector spaces, then $V$
becomes a $\gg$-module for the action:
\[ X\cdot v\equiv [\sigma(X),v],\qquad X\in\gg, v\in V.\]
Here we identify $V$ with its image in $\hh$. This action does not
depend on the splitting, and we can define an element $\omega\in
C^2(\gg;V)$ by:
\[ \omega(X,Y)=\sigma([X,Y])-[\sigma(X),\sigma(Y)],\qquad X,Y\in\gg.\]
Notice that the left-hand side actually lies in $\Ker\pi=V$. An
easy computation shows that $d\omega=0$, so we have a well-defined
class $[\omega]\in H^2(\gg;V)$. Again, it is not hard to check that
this class is independent of the splitting and that it vanishes iff
there exists a splitting of the sequence (as Lie algebras).

Now, the Second Whitehead Lemma states that if $\gg$ is a semisimple
Lie algebra and $V$ is any \emph{finite dimensional} $\gg$-module,
then $H^2(\gg;V)=0$. Certainly, we cannot apply this directly to the
short exact sequence (\ref{eq:seq:Lie:algebra}), since the the Lie
algebra $\M^2$ is neither finite dimensional nor abelian.  Still, just
like in the case of Lie algebra actions, we have:

\begin{theorem}[Weinstein \cite{Wein}]
If $\gg$ is a semisimple Lie algebra, then the Poisson bracket is formally
linearizable, i.e., there exists coordinates
$(x^1_\infty,\dots,x^m_\infty)$ around $x_0$ such that
\[  \{x^i_\infty,x^j_\infty\}=c^{ij}_k x^k_\infty+O(\infty),
\qquad (i,j=1,\dots, m).\]
\end{theorem}

For the proof, one looks again at the spaces $\M/\M^{k+1}$ of germs mod
terms of order higher than $k$. Let us look once more at the diagram
(but now of Lie algebras, rather than modules):
\[
\xymatrix{
\M\ar[r]\ar@<-1ex>[ddrrrrrr]_{\pi}&
\cdots\ar[r]&\M/\M^{k+1}\ar[r]^{\pi_k}&\M/\M^k\ar[r]&
\cdots\ar[r]&\M/\M^3\ar[r]^{\pi_2}&\M/\M^2\\
\\
&&&&&&T^*_{x_0}M=\gg_{x_0}\ar[uu]_{\phi_2}
\ar@{-->}[uul]_{\phi_3}\ar@{-->}[uulll]_{\phi_k}
\ar@{-->}[uullll]_{\phi_{k+1}}}
\]

We claim that one can construct the injective homomorphisms of Lie algebras
$\phi_k:\gg_{x_0}\to\M/\M^k$ by induction.

The case $k=2$ is trivial, so assume we have constructed $\phi_k$. Then,
we restrict the short exact sequence of (finite dimensional) Lie algebras
\[
\xymatrix{
0\ar[r]& \M^k/\M^{k+1}\ar[r]&\M/\M^{k+1}\ar[r]^{\pi_k}&\M/\M^k\ar[r]&0
},
\]
to the image $\phi_k(\gg)$. We obtain a new short exact sequence of
Lie algebras
\[
\xymatrix{
0\ar[r]& \M^k/\M^{k+1}\ar[r]&\pi_k^{-1}(\phi_k(\gg))/\M^{k+1}\ar[r]
&\phi_k(\gg)\ar[r]&0
},
\]
with $\M^k/\M^{k+1}$ abelian. Since $\gg$, and hence $\phi_k(\gg)$, is
semisimple, by the Second Whitehead Lemma there exists a splitting
$\sigma:\phi_k(\gg)\to \pi_k^{-1}(\phi_k(\gg))/\M^{k+1}$. Composing
with $\phi_k$ we obtain $\phi_{k+1}$.

By passing to the projective limit, we obtain a homomorphism of Lie
algebras
\[ \sigma: \gg_{x_0}\to \varprojlim \M/\M^k,\]
and, as above, this defines the coordinate system
$(x^1_\infty,\dots,x^m_\infty)$ that formally linearizes the Poisson
bracket.

\subsection{Analytic Linearization}
The proof given above for the formal linearization shows that
one constructs the formal linearizing coordinates by a succession of changes
of coordinates, which at each step removes terms of order $k$, without
modifying the terms of order less than $k$. Hence, it is natural to
ask if one can make the choices at each step so that, in the limit, one
obtains a smooth or analytic coordinate system, that linearizes the
action. In the analytic category this is indeed the case: 

\begin{theorem}[Conn \cite{Conn1}]
\labell{thm:poisson:conn1}
Let $\{\cdot,\cdot\}$ be an analytic Poisson structure which vanishes
at $x_0$. If the isotropy Lie algebra $\gg_{x_0}$ at $x_0$ is semisimple,
then there exist a local analytic coordinate system
$(x^1_\infty,\cdots,x^n_\infty)$ around $x_0$ in which the Poisson
structure is linear:
\begin{equation}
\{x^i_\infty,x^j_\infty\}=c^{ij}_k x^k_\infty.
\end{equation}
\end{theorem}

This is the analogue of the Guillemin and Sternberg linearization
result for Lie algebra actions. Again, we will sketch analytic and
geometric proofs. To simplify, we will denote by $\gg$ the isotropy
Lie algebra $\gg_{x_0}$.

\begin{proof}[Analytic proof] 
The analytic proof, due to Conn, uses a fast convergence method due to
Kolmogorov. To sketch a proof along these lines, let us look closer at
the formal linearization. We observe that the bracket relation
\[ [\M^k,\M^k]\subset \M^{2k-1},\]
implies that one can construct the sequence of linearizing coordinates 
so that, at each step, we remove terms of order $2^k$. 

To see this, let $\{X_1,\dots,X_n\}$ be a basis for the Lie algebra
$\gg$. We assume that we have constructed an injective homomorphism
$\phi_{\nu}:\gg\to \M/\M^{2^{\nu}}$ (\footnote{Note the change in
the indices!}), so that we have coordinates $(x^1_\nu,\cdots,x^n_\nu)$ 
where each $x^i_\nu$ represents the element $\phi_{\nu}(X_i)$ and,
moreover, the Poisson bracket satisfies:
\[
\{x^i_\nu,x^j_\nu\}= c^{ij}_k x^k_\nu+O(|x_\nu|^{2^{\nu}}).
\]
We want to construct the next iteration
$(x^1_{\nu+1},\cdots,x^n_{\nu+1})$ by the method we described above. 

For that, we note that $\gg$ acts on on the spaces
$\M^{2^\nu}/\M^{2^{\nu+1}}$ by Hamiltonian diffeomorphisms:
\begin{equation}
X_i\cdot f
=c^{ij}_k x^k_\nu\frac{\partial f}{\partial x^j_\nu} 
=\{x^i_{\nu},f\}
\quad {\rm mod} \,\M^{2^{\nu+1}}.
\end{equation}
We consider the remainder term $R_\nu:\wedge^2 \gg \to
\M^{2^\nu}/\M^{2^{\nu+1}}$ which is given by:
\[ R_\nu(X_i\wedge X_j)=\{x^i_\nu,x^j_\nu\}-c^{ij}_k
x^k_\nu \quad {\rm mod}\, \M^{2^{\nu+1}}.\] 
The Jacobi identity shows that $R_\nu$ is a cocycle and so defines a 
cohomology class in 
$H^2(\gg;\M^{2^\nu}/\M^{2^{\nu+1}})$. Since $\gg$ is
semisimple, by the second Whitehead lemma, this class must vanish. 
Let us denote by $\sigma_{\nu+1}:\gg\to \M^{2^\nu}/\M^{2^{\nu+1}}$ a
primitive. We can choose local analytic functions $f^i_{\nu+1}$
which represent the elements $\sigma_{\nu+1}(X_i)\in
\M^{2^\nu}/\M^{2^{\nu+1}}$. The new coordinates are then defined by
\begin{equation}
x^i_{\nu+1}\equiv x^i_\nu-f^i_{\nu+1}\,.
\end{equation}
and it is easy to verify that:
\[
\{x^i_{\nu+1},x^j_{\nu+1}\}=c^{ij}_k x^k_{\nu+1}+O(|x_\nu|^{2^{\nu+1}}).
\]
This establishes the claim.

Now we need to look at the formal limit of these systems of
coordinates, which we denote by $(x^1_\infty,\cdots,x^n_\infty)$. 
To complete the proof of Theorem \ref{thm:poisson:conn1}, we ``just''
have to show that these formal coordinates are local
analytic coordinates.

In order to make estimates in the algorithm above, we need a normed
version of the second Whitehead lemma. For this, we denote by
$\gg_\Cc$ the complexification of $\gg$. This complex Lie algebra is
semisimple, and we denote by $\gg_0$ its compact real form.  We
consider a (finite dimensional) Hermitian complex vector space which
is a Hermitian $\gg_0$-module (i.e., $\gg_0$ acts by Hermitian
transformations), and we let $\|~\|$ denote the norm corresponding 
to the Hermitian structure. Then:

\begin{proposition}[Conn \cite{Conn1}]
\labell{prop:vanishing:analytic}
There exists a linear map 
$h:\wedge^2\gg_\Cc\otimes W\to\gg_\Cc\otimes W$ 
and a positive constant $D$ (which depends only on $\gg$) such that, for
every $R\in\wedge^2\gg_\Cc \otimes W$,
\[
\|h(R)\|\leq D\|R\|,
\]
and
\[
d(R)=0\quad \Longrightarrow\quad R=d\big( h(R)\big).
\]
where $d$ is the differential of the complex $(C^\bullet(\gg_0;W),d)$. 
The homotopy operator $h$ is real, i.e., maps real cochains to real cochains.
\end{proposition}

We would like to apply this proposition to the spaces
$\M^{2^\nu}/\M^{2^{\nu+1}}$, so we must furnish them
with Hermitian metrics which are preserved by the $\gg_0$-action. This can be 
done as follows:

The Lie algebra $\gg$ acts on $T_{x_0}M=\gg^\ast$ by the coadjoint
action. We can complexify this action obtaining an action of
$\gg_\Cc$: we choose our original local coordinates so that their
complexification $(z_1,\dots,z_n)$ determined a basis for $\gg_\Cc^*$
which is orthonormal with respect to the opposite of the Killing form
of $\gg_0$ (which is a Hermitian metric of $\gg_\Cc^\ast$). This
Hermitian metric is preserved by the coadjoint action of the compact 
Lie algebra $\gg_0$.

Let $r$ be a positive number. For each integer $\nu$, for which we
have constructed analytic coordinates $(z^1_\nu,\cdots,z^n_\nu)$, 
we define the (deformed) balls
\[
B_{\nu,r}=\{w\in\Cc^n\,|\, \sqrt{\sum_i|z^i_\nu(w)|^2}\leq r\}
\]
In this ball we take the standard volume form $d\mu_\nu$ relative
to the coordinates $(z^1_\nu,\cdots,z^n_\nu)$ and we denote
by $V_r$ the volume of the ball. This ball is analytically
diffeomorphic to the standard ball of radius $r$ via the coordinate
system $(z^1_\nu,\cdots,z^n_\nu)$, and we
define the $L^2$-norm on local analytic functions:
\[
\|f\|_{\nu,r}=\sqrt{\frac{1}{V_r}\int_{B_{\nu,r}} |f(w)|^2
d\mu_\nu}
\]
If we identify the quotients $\M^{2^\nu}/\M^{2^{\nu+1}}$ with space of
polynomials functions, these norms determined corresponding norms on
the quotient. One can show that the norms obtained in this way are
Hermitian. For example, the norm $\|f\|_{0,r}$ is given by the
following classic Hermitian metric: if $f=\sum_{\alpha\in\Nn^n}
a_\alpha \prod (z^i)^{\alpha_i}$ and $g=\sum_{\alpha\in\Nn^n} b_\alpha
\prod (z^i)^{\alpha_i}$, then we have
\[
\langle f,g\rangle_{0,r}=\sum_{\alpha\in\Nn^n}
\frac{\alpha!n!}{(|\alpha|+n)!} a_\alpha \overline{b}_\alpha
r^{2|\alpha|},
\]
where $\alpha!=\prod_i\alpha_i!$ and $|\alpha|=\sum\alpha_i$.
The others norms are also given by similar Hermitian metrics.

Therefore, the spaces $\M^{2^\nu}/\M^{2^{\nu+1}}$
equipped with the norms $\|\,\|_{\nu,r}$ are Hermitian spaces.
Now remark that, since the action of $\gg_0$ on $\Cc^n$ preserves
the Hermitian metric of $\Cc^n$, its action on
$\M^{2^\nu}/\M^{2^{\nu+1}}$ preserves the norm
$\|~\|_{\nu,r}$. We can then apply the proposition above to conclude
that:
\[ 
\max_{i}\|z^i_{\nu+1}-z^i_{\nu}\|_{\nu,r}\le D 
\max_{i,j}\|\{z^i_\nu,z^j_\nu\}-c^{ij}_k z^k_\nu\|_{\nu,r}
\]

Using this estimates, one constructs a decreasing sequence of radii 
${(r_\nu)}_\nu$ such that the balls $B_{\nu,r_\nu}$ are well-defined
and there exists an integer $K$ and open neighborhood $U$ of the
origin, such that: 
\[ U\subset \bigcap_{\nu\geq K} B_{\nu,r_\nu}.\]
Moreover, for every $w\in U$ and $\nu\geq K$, we have uniform bounds
on the norms of the corrections:
\[ \sup_{w\in B_{\nu,r_\nu}}|z^i_\nu(w)-z^i_{\nu+1}(w)|,\]
and on the Poisson structure. We conclude that the limit coordinate system
$(z^1_\infty,\cdots,z^n_\infty)$ is analytic in $U$, and therefore
furnishes analytic linearizing coordinates.
\end{proof}

\begin{proof}[Geometric proof]
In \cite{Zung2}, Zung sketches a geometric proof of the theorem. This
proof consists of the following steps:
\begin{enumerate}
\item[(i)] Construct a basis of Casimirs functions
  $\{\tilde{C}_1,\dots,\tilde{C}_r\}$ of the Poisson bracket $\{~,~\}$
  by perturbing a basis of Casimir functions $\{{C}_1,\dots,{C}_r\}$
  of the linear bracket. This is possible because the Casimirs of the
  linear bracket are given as certain invariant polynomials of
  integrals of the leafwise symplectic form over the periods
  (2-cycles) of the coadjoint orbits. By analyticity, one can perturb
  these 2-cycles to 2-cycles of the symplectic leaves of the
  non-linear bracket.
\item[(ii)] Use the Casimir functions, to conclude that the
  symplectic foliation of the original Poisson bracket is analytically
  diffeomorphic to the symplectic foliation of the linear
  bracket. Hence, one can assume that they are diffeomorphic.
\item[(iii)] Use averaging and Moser's path method to deform the
  Poisson bracket to the linear bracket along an isotopy which is the
  time-1 flow of a certain time-dependent vector field.
\end{enumerate}
Not every detail is given in \cite{Zung2}, but it is likely that a
complete proof can be given along these lines.
\end{proof}

\subsection{Smooth Linearization}
Let us now turn to the linearization problem of smooth Poisson
brackets. In the smooth setting it is not enough to assume that the
$\gg_{x_0}$ is semisimple as one sees from the following example, due to
Weinstein (see \cite{Wein}).

\begin{example}
Consider the Lie algebra $\mathfrak{sl}(2)$ and identify its dual with
$\Rr^3$, with coordinates $\{x,y,z\}$. Then we have the linear Poisson
bracket: 
\[ \{x,y\}=-z,\quad \{y,z\}=x, \quad \{z,x\}=y.\]
The leaves of this Poisson bracket are the level sets of the quadratic form
$x^2+y^2-z^2=r^2-z^2$. On the other hand, we can perturb this linear
bracket to the non-linear bracket:
\begin{align*} 
\{x,y\}&=-z,\\
\{y,z\}&=x-\frac{y}{r^2}g(r^2-z^2),\\
\{z,x\}&=y+\frac{x}{r^2}g(r^2-z^2),
\end{align*}
where $g\in C^\infty(\Rr)$ is such that $g(x)>0$, if $x>0$, and $g(x)=0$,
if $x\le 0$. Now, observe that the hamiltonian vector field for the
function $h(x,y,z)=z$ is given by:
\[ 
X_h=x\frac{\partial}{\partial y}-y\frac{\partial}{\partial x}
-\frac{g(r^2-z^2)}{r^2}
\left(x\frac{\partial}{\partial x}
+y\frac{\partial}{\partial y}\right),\]
and has integral curves outside the cone $r=z$ which spiral towards
the cone. Since no hamiltonian vector field in $\mathfrak{sl}(2)^*$
has this behavior, this Poisson bracket is not linearizable.
\end{example}

Again, like the Lie algebra case, we have:

\begin{theorem}[Conn \cite{Conn2}]
Let $(M,\{~,~\})$ be a smooth Poisson manifold vanishing at $x_0$. If
the isotropy Lie algebra $\gg_{x_0}$ is semisimple of compact type, then
there exists smooth local coordinates $(x^1_\infty,\dots,x^m_\infty)$ around
$x_0$ such that
\[  \{x^i_\infty,x^j_\infty\}=c^{ij}_k x^k_\infty \qquad (i,j=1,\dots, m).\]
\end{theorem}

There is also a $C^p$ version of this linearization theorem (see
\cite{MoZung}). Again, we shall sketch an analytic and a geometric
proof of this result.

\begin{proof}[Analytic proof] As in the analytic case, we would like
to construct a sequence of smooth local coordinates systems
$(x^1_\nu,\cdots,x^n_\nu)$ which converges to a smooth local
coordinates system $(x^1_\infty,\cdots,x^n_\infty)$ in which the
Poisson structure is linear. In the proof of the analytic case, at
each step, we used a homotopy operator to correct the remainder terms:
\[ R_\nu:\wedge^2 \gg \to \M^{2^\nu}/\M^{2^{\nu+1}}.\]
We could do this because the $R_{\nu}$ were 2-cocycles in the
Chevalley-Eilenberg complex
$C^\bullet(\gg;\M^{2^\nu}/\M^{2^{\nu+1}})$. Here we see the first
major difference from the analytic case: the bracket maybe non-linear, while
all such remainder terms vanish! Hence, we are now forced to work on the
$\gg$-modules $\M^{2^\nu}$, and put aside the quotients
$\M^{2^\nu}/\M^{2^{\nu+1}}$. 

So let us consider now the remainder terms as maps:
\[ R_\nu:\wedge^2 \gg \to \M^{2^\nu}.\]
These are not cocycles anymore. However, we will think of them as
\emph{almost cocycles} and we apply the homotopy operator, so we now
must control the size of the resulting \emph{almost coboundary} $h(R_\nu)$.
Here again, we use a normed version of the Whitehead lemma for a
certain orthogonal $\gg$-module (of infinite dimension).

So we denote by $(x^1,\cdots,x^n)$ a coordinate system centered at
$x_0$, such that $\{d_{x_0}x^1,\cdots,d_{x_0}x^n\}$ is an orthonormal
basis of $\gg$, with respect to the negative of the Killing form
(since $\gg$ is compact, this is positive definite). We have
\[ [d_{x_0}x^i,d_{x_0}x^j]=c^{ij}_k d_{x_0}x^k.\] 
Let $r$ be a positive real number and denote by $B_r$ the closed ball
of radius $r$ and centered at 0, with respect to the
coordinates $(x^1,\cdots,x^n)$. We denote by $\mathcal{C}_r$ the
space of smooth functions on $B_r$, vanishing at 0 and whose first
derivatives also vanish at 0. The Lie algebra $\gg$ acts on $B_r$
via the coadjoint action, and this induces an action on
$\mathcal{C}_r$. On the space $\mathcal{C}_r$ the
Sobolev $H_k$-metrics are defined by:
\[
\langle f,g \rangle_{k,r} :=\sum_{|\alpha|\leq k} \int_{B_r}
\big(\frac{|\alpha|!}{\alpha!}\big)
\Big(\frac{\partial^{|\alpha|}f}{\partial x^\alpha}(z)\Big) {\rm
d}\mu(z)
\]
and we denote by $\|\,,\,\|_{k,r}$ the corresponding norm. The action
of $\gg$ on $\mathcal{C}_r$ preserves these norms. The $\gg$-module
$(\mathcal{C}_r,\langle\,,\, \rangle_{k,r})$ is of infinite
dimension, but Conn showed that the complex associated to this
$\gg$-module shares the same properties as the finite
dimensional $\gg$-modules, that we have consider in Proposition
\ref{prop:vanishing:analytic}. More precisely, there exists a 
homotopy operator $h$ for the truncated Chevalley-Eilenberg complex:
\[ \mathcal{C}_r \otimes \wedge^1 \gg^\ast \stackrel{h}{\leftarrow} 
\mathcal{C}_r \otimes \wedge^2 \gg^\ast \stackrel{h}{\leftarrow} 
\mathcal{C}_r \otimes
\wedge^3 \gg^\ast \]
such that
\[
d \circ h + h \circ d = {\rm Id}_{\mathcal{C}_r
 \otimes \wedge^2\gg^\ast},
\]
and a positive constant $C$, which is independent of
the radius $r$, such that
\begin{equation}
\|h(u)\|_{k,r}^H \leq C \|u\|_{k,r}^H,\quad \forall k \geq 0.
\label{eqn:estimate} 
\end{equation}
(the same property holds for every infinite direct sum of
finite-dimensional orthogonal $\gg$-modules).

The homotopy operator is well-controlled by Sobolev norms. However,
in order to control the differentiability at each step, we have to use the
$C^k$-norms : 
$$\|F\|_{k,r}:=\sup_{|\alpha|\leq k}\sup_{z\in B_r} |D^\alpha F(z)|\,.$$
So applying Sobolev's lemma, the inequality (\ref{eqn:estimate}) leads
to the estimate:
\begin{equation}
\|h(u)\|_{k,r} \leq M \|u\|_{k+s,r} \, 
\label{eqn:estimate2} 
\end{equation}
where $M$ is a positive constant and $s$ a fix positive integer. So
now we are facing a loss of differentiability: to control the
$C^k$-differentiability of $h(u)$, we need to control the
$C^{k+s}$-differentiability of $u$. 

There is a standard way to get around difficulty using \emph{smoothing
operators} on $\mathcal{C}_{r}$:
One can construct, a family of linear operators $S(t)=S_r(t)$ $(t>1)$ from
$\mathcal{C}_{r}$ into itself, satisfying the following properties:
\begin{enumerate}
\item[(i)] $\|S(t)f\|_{p,r}\leq C_{p,q} t^{(p-q)}\| f \|_{q,r}$ and,
\item[(ii)] $\|(I-S(t))f\|_{q,r}\leq C_{p,q} t^{(q-p)}\| f\|_{p,r}$.
\end{enumerate}
for any $f \in \mathcal{C}_r$. Here $p,q$ are any non-negative
integers such that $p \geq q$, $I$ denotes the identity map, and
$C_{p,q}$ denotes a constant which depends on $p$ and $q$. Note that
(ii) means that $S(t)$ is close to identity and converges
to the identity when $t\to\infty$. On the other hand, (i) means
that $f$ becomes ``smoother'' when we apply $S(t)$ to it. 

We are now ready to proceed with the iteration. Assume we have
constructed coordinates $(x^1_\nu,\cdots,x^n_\nu)$, so that 
the Poisson bracket satisfies:
\[
\{x^i_\nu,x^j_\nu\}= c^{ij}_k x^k_\nu+O(|x_\nu|^{2^{\nu}}).
\]
We want to construct the next iteration
$(x^1_{\nu+1},\cdots,x^n_{\nu+1})$.
We consider the remainder term $R_\nu:\wedge^2 \gg \to \M^{2^\nu}$ 
which is given by:
\[ R_\nu(X_i\wedge X_j)=\{x^i_\nu,x^j_\nu\}-c^{ij}_k
x^k_\nu.\] 
Then $R_{\nu}$ is an almost cocycle: $dR_\nu$ is a quadratic function in
$R_\nu$ i.e., if $R_\nu$ is ``$\varepsilon$-small'' then $dR_\nu$ is
``$\varepsilon^2$-small''. We define the next change of coordinates to
be given by the diffeomorphism $\phi_{\nu+1}$ defined by
\begin{equation}
\phi_{\nu+1}={\rm Id}-S(t_\nu)(h(R_\nu),
\end{equation}
where the real numbers $t_\nu$ are defined by 
$t_{\nu+1}=t_\nu^{3/2}$ with $t_0>1$ (this choice of smoothing 
parameter is standard). Now we are left with the question of
convergence, and here arrives the most technical part. Ones constructs
a decreasing sequence of  radii $\{r_\nu\}$, which
converges to $r>0$, such that the corrected Poisson bracket
$\{\,,\,\}_\nu$ is well defined on $B_{r_\nu}$, and such that the
$\|R_\nu\|_{k,r_\nu}$ and the $\|\phi_\nu-Id\|_{k,r_\nu}$ converge
to 0 exponentially fast. This is proved using some hard estimates that
follow from (\ref{eqn:estimate2}) and the properties of the smoothing
operators.
\end{proof}

\begin{remark}
This method of the proof is known as the ``Nash-Moser method''. A
sketch of a more conceptual proof, using the Nash-Moser Inverse
Function Theorem on tame Fr\'echet spaces (see
\cite{Hamilton-NashMoser1982}) is given by Desolneux-Moulis in
\cite{Desou}. However, this proof seems to be incomplete.
\end{remark}

\begin{proof}[Geometric proof]
In \cite{CrFe2}, Crainic and Fernandes propose a soft geometric 
proof of the theorem, along the lines of the  infinitesimal version of
Bochner's Theorem (Theorem \ref{thm:infin:Bochner}). The proof has two
major steps: 
\begin{enumerate}
\item[(i)] Integrate (a neighborhood) of the Poisson manifold to a
  symplectic groupoid.
\item[(ii)] The resulting symplectic groupoid has compact isotropy
  group $G_{x_0}$ so we can apply Theorem \ref{thm:proper:groupoid} to
  linearize it.
\end{enumerate}
We refer to \cite{CrFe2} for more details. Also, in \cite{CrFeMo} this
result is interpreted in the context of deformation cohomology.
\end{proof}

\begin{remark}
It is natural to try and use the approach of Flato, Pinczon and Simon
(\cite{Flato}) mentioned in Remark \ref{rem:flato}, to simplify Conn's theorems. 
Unfortunately, so far, no such simplification was obtained (see the remarks in
\cite{Wein:linearization}).
\end{remark}

As we have mentioned in the introduction, an important open problem is
to characterize which Lie algebras are non-degenerate relative to
Poisson structures.  In dimension 2, there are two Lie algebras (up to
isomorphism): the 2-dimensional abelian Lie algebra (as any abelian
Lie algebra) is degenerate, while Arnold has shown in \cite{Arnold}
that the non-abelian 2-dimensional Lie algebra is non-degenerate. The
case of dimension 3 is solved completely by Dufour in \cite{Dufour:Poisson},
while a discussion of dimension 4 is given by Molinier in
\cite{Mol}. However, in general, this problem is completely open.

\begin{remark}
In the formal case, the (non)degeneracy of a Lie algebra should be
related to its (non)rigidity (recall that $H^2(\gg;\gg)=0$ implies
$\gg$ rigid). in \cite{Bordemann}, Bordemann, Makhlouf and Petit show 
that, if the universal enveloping algebra $\mathcal{U}(\ll)$
is rigid as an associative algebra, then $\ll$ is formally
non-degenerate.
\end{remark}

\section{Levi decomposition}
\labell{sec:Levi}

A Poisson bracket defines a Lie algebra structure on
$C^\infty(M)$. This infinite dimensional Lie algebra
retains a finite dimensional flavor, due to the Leibniz identity:
\[ \{f,gh\}=\{f,g\}h+g\{f,h\}.\]
Therefore, it is natural to seek a Levi decomposition for this Lie
algebra. We study in this section this decomposition, which may also
be seen as a semi-linearization generalizing the linearization results
of the previous section.

As we saw in the previous sections, the linearization of Lie algebra
actions involves the first Whitehead lemma, whereas the linearization
of Poisson structures involves the second Whitehead lemma. We will see
that the Levi decomposition of Poisson structures involves both of
them so, somehow, it combines both the linearization of Lie algebra
actions and the linearization of Poisson structures.

Let us first recall the Levi decomposition for finite dimensional Lie
algebras. Let $\gg$ be a finite-dimensional Lie algebra and denote by
$\rr$ its radical (i.e., its maximal solvable radical). The quotient
Lie algebra $\gg/\rr$ is then semisimple and we can write the
following exact sequence:
\[ \xymatrix{ 0\ar[r]& \rr\ar[r]& \gg\ar[r]&
      \gg/\rr\ar[r] \ar[r]&0}.\]
The Levi-Malcev theorem states that this exact sequence of Lie
algebras admits a splitting $\sigma:\gg/\rr\to\ll$. If
we denote by $\ss$ the image $\sigma(\gg/\rr)$, then we can write
$\gg$ as a semi-direct product:
\[\gg= \ss \ltimes \rr,\]
which is called the \textbf{Levi decomposition} of $\gg$.  In general,
for infinite dimensional Lie algebras, there exists no Levi-Malcev
decomposition. However, the Levi-Malcev decomposition does hold for
filtered, pro-finite, Lie algebras (details are given in the upcoming
book \cite{Book:JP-Zung}).

Let us turn now to the case of a Poisson manifold. Let $\{~,~\}$ be
a Poisson bracket on a manifold $M$ vanishing at $x_0\in M$. Denote by 
$\gg=\gg_{x_0}=T^*_{x_0}M$ the isotropy Lie algebra at $x_0$. Now we
choose $\ss\subset\gg$ to be:
\begin{itemize}
\item a Levi factor of $\gg$, in the formal or analytic case, or
\item a maximal compact semisimple Lie subalgebra of $\gg$ (a
  \emph{compact Levi factor}) in the smooth case.
\end{itemize}
As before, we denote by $\M$ the maximal ideal formed by functions that vanish
at $x_0$. The Poisson bracket defines a Lie algebra structure on $\M$
and we let $\R\subset \M$ be defined by:
\[ \R=\{f\in \M:d_{x_0}f\in\ss\}.\]
We should think of $\R$ as the \textbf{radical} of the Lie algebra
$\M$. Now we have a short exact sequence of Lie algebras:
\begin{equation}
    \labell{eq:seq:Lie:brackets}
    \xymatrix{ 0\ar[r]& \R\ar[r]& \M\ar[r]&
      \ss\ar[r] \ar@(dl,dr)@{-->}[l]^{\sigma}&0}
\end{equation}
We shall call a  splitting $\sigma$ of this short exact sequence of Lie
algebras a \textbf{Levi decomposition for the Poisson bracket}. 

Notice that when $\gg=\ss$ we obtain the standard linearization
problem for Poisson brackets that we have studied in the previous section.

Let us give the coordinate version of the Levi decomposition for
Poisson brackets. Under the assumptions above we can choose a
complement $\rr$ to $\ss$ in $\gg$, such that
\[ [\gg,\rr]\subset \rr.\]
Therefore, we can choose coordinates $(x^1,\cdots,x^m,y^1,\cdots,y^{n-m})$ 
around $x_0$ such that $(d_{x_0}x^1,\cdots,d_{x_0}x^m)$ spans $\ss$
and $(d_{x_0}y^1,\cdots,d_{x_0}y^{n-m})$ spans $\rr$. In terms of the
Poisson bracket this means that:
\begin{align*}
\{x^i,x^j\}&= c^{ij}_k x^k + O(2),\\
\{x^i,y^j\}&= a^{ij}_k y^k + O(2),\\
\{y^i,y^j\}&= O(1).
\end{align*}
Now the existence of a (formal, analytic, smooth) Levi decomposition
for the Poisson brackets amounts to the existence of a (formal,
analytic, smooth) coordinate system 
$(x^1_\infty,\cdots,x^m_\infty,y^1_\infty,\cdots,y^{n-m}_\infty)$, 
for which the higher order terms in the first and second pair of
brackets vanish:
\begin{align*}
\{x^i_\infty,x^j_\infty\}&= c^{ij}_k x^k_\infty,\\
\{x^i_\infty,y^j_\infty\}&= a^{ij}_k y^k_\infty,
\end{align*}
(note that nothing is said about the third pair of brackets)

Now, exactly as before, we have the following results:
\begin{theorem}[Wade \cite{Wade}]
Any Poisson bracket that vanishes at $x_0$ admits a formal Levi decomposition.
\end{theorem}

\begin{theorem}[Zung \cite{Zung}]
Any analytic Poisson bracket that vanishes at $x_0$ admits an analytic
Levi decomposition
\end{theorem}

\begin{theorem}[Monnier and Zung \cite{MoZung}]
Any smooth Poisson structure that vanishes at $x_0$ admits a smooth
Levi decomposition. 
\end{theorem}

If in the above theorems the full isotropy Lie algebra $\gg$ 
is semisimple (in the formal or analytic case) or compact semisimple
(in the smooth case), we rediscover the linearization theorems
given in the previous section. Therefore, it should be no surprise
that the proofs of the Levi decompositions are higher order versions
of the proofs of the linearization results for Poisson brackets.

For example, for the analytic proofs, one constructs a sequence of
(formal, analytic, smooth) coordinates  systems
$(x^1_\nu,\cdots,x^m_\nu,y^1_\nu,\cdots,y^{n-m}_\nu)$ 
which converges (formal, analytic, smooth) to the coordinate system 
$(x^1_\infty,\cdots,x^m_\infty,y^1_\infty,\cdots,y^{n-m}_\infty)$.
The additional complication is that, at each step, we have to deal
with the $\{x^i,x^j\}$-terms and the $\{x^i,y^\alpha\}$-terms
together. The algorithm to linearize the $\{x^i,x^j\}$-terms is
the similar to the one described in Section \ref{sec:Poisson:brackets},
and involves the vanishing of the second Lie algebra cohomology of $\ss$. To
linearize the $\{x^i,y^\alpha\}$-terms, one uses a procedure similar
to the one use in Section \ref{sec:Lie:algebras}, and which involves,
now, the vanishing of the first Lie algebra cohomology of $\ss$. In
the analytic case or the smooth case, the estimations get more
involved.

\begin{remark}
At least in the formal case, instead of treating both
$\{x^i,x^j\}$-terms and the $\{x^i,y^\alpha\}$-terms at each step, one
can first linearize the $\{x^i,x^j\}$-terms. This leads to an action of
$\ss$ on $\M$ which can then be linearized, thus taking care of the
$\{x^i,y^\alpha\}$-terms.
\end{remark}

\begin{example}
As we pointed out above, the Levi decomposition can be thought of as a
semi-linearization and may be seen as a first step in the tentative of
linearizing Poisson structures whose linear part does not satisfy
the assumptions of the linearization results of Section
\ref{sec:Poisson:brackets}. 

For example, let $\gg=\mathfrak{aff}(n)$ be the Lie
algebra of affine transformations of $\Rr^n$. This Lie algebra has
radical $\rr=\Rr\ltimes\Rr^n$, where $1$ acts on $\Rr^n$ by the
identity map, and it admits the Levi factor
$\ss=\mathfrak{sl}(n)$. Given an analytic Poisson bracket $\Pi$ with a
fixed point $x_0$ and isotropy $\gg_{x_0}=\mathfrak{aff}(n)$, we can apply the
Levi decomposition to construct semi-linearizing coordinates
$x_1,\dots,x_{n^2-1},y_0,\dots,y_n$. In \cite{JP-Zung}, it is shown
that one can choose the coordinates so that one also has:
\[ \{y_0,x_i\}=0,\qquad \{y_0,y_j\}=y_j.\]
Denoting by $F_1(x),\dots,F_n(x)$ the elementary symmetric polynomials
in $\mathfrak{gl}(n)$, with a bit more work, one can use these
relations to show that there exist analytic functions $f_i=f_i(x)$
so that the vector field 
\[ Y=\sum_{i=1}^n f_i X_{F_i},\]
satisfies:
\[ \Lie_Y(\Pi)=\Pi_0-\Pi.\]
where $\Pi_0$ is the linear part of $\Pi$. This is explain in detail
in \cite{JP-Zung}.

Now consider the path of Poisson structures:
\[ \Pi_t(x,y)\equiv\frac{1}{t}\Pi(t(x,y)),\]
which joins $\Pi$ to $\Pi_0$. The time-dependent vector field
$Y_t(x,y)=\frac{1}{t^2}Y(t(x,y))$ satisfies 
\[ \Lie_{Y_t}\Pi_t=\frac{d\Pi_t}{dt},\]
so we conclude that the time-1 flow of $Y_t$ is a Poisson
diffeomorphism between $\Pi_0$ and $\Pi$ defined in a neighborhood of $x_0$.
In other words, $\mathfrak{aff}(n)$ is analytically non-degenerate.
\end{example}

\section{Linearization and Levi decomposition of Lie algebroids}
\labell{sec:algebroids}

The similarity between the linearization of Lie algebra actions and
Lie algebroids, that we have seen before, can be fully understood in the
context of Lie algebroids. This was first observed by Weinstein in
\cite{Wein:linearization}. In this section, we will consider only the
smooth case. The formal and the analytic case are similar.

Let $\pi:A\to M$ be a Lie algebroid with Lie bracket $[~,~]$ on
$\Gamma(A)$ and anchor $\#:A\to TM$. For background on Lie algebroids
we refer to the book by Cannas da Silva and Weinstein \cite{CaWe}. The
splitting theorem for Lie algebroids (see
\cite{Dufour:algebroids,Fer}) states that, for every $x_0\in M$,
there exists a neighborhood $U$, such that $A|_U$ is a direct sum of a
Lie algebroid with constant rank anchor and a Lie algebroid whose
anchor vanishes at $x_0$.  So, henceforth, we assume that $x_0\in M$
is a fixed point, i.e., $\#_{x_0}\equiv 0$.

If we choose a basis $\{e_1,\cdots,e_r\}$ of local sections
of $A$ and a system of coordinates $(x^1,\cdots,x^n)$ around
$x_0$, we have:
\begin{align*}
[e_i,e_j]&=(c^k_{ij}+O(1)) e_k, \\
\# e_i&=(b_{ij}^k x^j+ O(2))\frac{\partial}{\partial x^k}.
\end{align*}
where $b_{ij}^k$ and $c^k_{ij}$ are certain structure constants. So we
have now the following \textbf{linearization problem for Lie algebroids}:
\begin{itemize}
\item Is there a choice of basis $\{e_1^\infty,\cdots,e_r^\infty\}$
  and a choice of local coordinates $(x^1_\infty,\cdots,x^n_\infty)$ where
  the higher order terms vanish?
\end{itemize}

The linear part of the Lie algebroid defined by
\begin{align*}
[e_i,e_j]&=c^k_{ij} e_k, \\
\# e_i&=b_{ij}^k x^j\frac{\partial}{\partial x^k}.
\end{align*}
is a special kind of Lie algebroid, called an \textbf{action Lie
algebroid}. In fact, the first set of relations defines a Lie algebra
$\gg_{x_0}$, called the isotropy Lie algebra at $x_0$, while the 
second set of relations defines a linear action of
$\gg_{x_0}$ on $T_{x_0}M$. This Lie algebra and this action are, in
fact, independent of the choice of basis and of the choice of
coordinates. We recall that \emph{any} Lie algebra action
$\rho:\gg\to\X(M)$ defines an action Lie algebroid: the bundle
$A=M\times\gg$ is the trivial vector bundle with fiber $\gg$, the
anchor $\#:A\to TM$ is given by $\#_x v=\rho(x)\cdot v$, and the
Lie bracket is defined fiber-wise on constant sections, and extended to
any sections so that it satisfies the Leibniz identity. The linear
part is an example of a linear action Lie algebroid, which we denote
by $A^L=\gg\ltimes T_{x_0}M$.

Note that if we eliminate the higher order terms in the Lie bracket,
then the second set of equations defines a (non-linear) Lie algebra
action of $\gg_{x_0}$. The remaining issue is to linearize the action,
which, under good assumptions (semisimplicity, compact type) on the
Lie algebra $\gg_{x_0}$, can be done using the results of Section
\ref{sec:Lie:algebras}.

Now, one has linearization results entirely analogous to what we have
seen in the Lie algebra and Poisson cases. For example, in the smooth
case we have:

\begin{theorem}[Monnier and Zung \cite{MoZung}]
Any smooth Lie algebroid with a fixed point $x_0$ for which the
isotropy $\gg_{x_0}$ is semisimple of compact type is smoothly
linearizable.
\end{theorem}

Similarly, if $\gg_{x_0}$ is semisimple, we have formal linearization
(Weinstein \cite{Wein:linearization}, Dufour \cite{Dufour:algebroids})
and analytic linearization (Zung \cite{Zung}).

\begin{example}
It is easy to see that, for any Lie algebra $\gg$, to linearize an
action $\rho:\gg\to\X(M)$ around a fixed point $x_0\in M$ is
equivalent to linearize the action Lie algebroid $\gg\ltimes_\rho TM$
around $x_0$.  Similarly, for a manifold $M$, to linearize a Poisson
bracket $\{~,~\}$ around a zero $x_0\in M$ one can show that is
equivalent to linearize the cotangent Lie algebroid $T^*M$. Therefore,
these results imply the linearization results of the previous
sections.
\end{example}

One can also look for a \textbf{Levi decomposition for Lie algebroids}. 
In the smooth case, we choose a decomposition
\[ \gg_{x_0}=\ss+\rr,\]
where  $\ss$ is a compact Levi factor and $\rr$ is a subspace
invariant under the adjoint action. 

\begin{theorem}[Monnier and Zung \cite{MoZung}]
Let $A$ be a smooth Lie algebroid with anchor $\#:A\to TM$ with
$\#_{x_0}=0$, and fix a decomposition of the isotropy
$\gg_{x_0}=\ss+\rr$, as above. Then, there exists a local basis of
sections $(e_\infty^1,\cdots,e_\infty^m,f_\infty^1,\cdots,f_\infty^{r-m})$
($m=\dim\ss$) and a local system of coordinates
$(x_\infty^1,\cdots,x_\infty^n)$ around $x_0$, such that:
\begin{align*}
[e^\infty_i,e^\infty_j]&= c_{ij}^k e_\infty^k, \\ 
[e^\infty_i,f^\infty_j]&= a_{ij}^k v_\infty^k, \\
\#e^\infty_i &= b_{ij}^k x_\infty^l\frac{\partial}{\partial x_\infty^k},
\end{align*}
where $c^{ij}_k,a^{ij}_k,b^{ij}_k$ are constants, with $c^{ij}_k$ being the
structural constants of the compact semisimple Lie algebra $\ss$.
\end{theorem}

Similarly, if $\ss$ is semisimple, we have formal and analytic Levi
decompositions (see \cite{Zung} for details).

The proofs of these results can be reduced to the proofs of the
corresponding statements for Poisson structures. Indeed, remind that
any Lie algebroid $A$ induces (and is, in fact, determined by) a
fiber-wise linear Poisson structure on the dual bundle $A^\ast$.  More
precisely, if $(x^1,\cdots,x^n)$ is a local coordinate system and
$(e_1,\cdots,e_r)$ is a local basis of sections, then we can think of
$(x^1,\cdots,x^n,e_1,\cdots,e_r)$ has a coordinate system for
$A^\ast$, which is linear on the fibers. The Poisson structure on
$A^*$ is given by
\begin{align*}
\{e_i,e_j\} & = [e_i,e_j],\\ 
\{e_i,x^j\} & = \#e_i(x^j),\\ 
\{x^i,x^j\} & = 0.
\end{align*}
This Poisson structure is fiber-wise linear in the following sense:
\begin{enumerate}
\item[(i)] the bracket of two fiber-wise linear functions is again a
  fiber-wise linear function,
\item[(ii)] the bracket of a fiber-wise linear function with a basis
  function is a basic function, and
\item[(iii)] the bracket of two basic functions is zero.
\end{enumerate}
The proof of the theorems above consist in showing that the Levi
decomposition for the fiber-wise linear Poisson structure on $A^*$
yields the Levi decomposition of the Lie algebroid $A$, around
$x_0$. The proof is similar to the Poisson case but one must modify
the $\gg$-modules in order to preserve the ``fiber-wise linear''
feature.

\bibliographystyle{amsplain}

\begin{thebibliography}{11}

\bibitem{Arnold} V.I.~Arnold, \emph{Geometrical methods in the theory
of ordinary differential equations}, Spinger-Verlag, New York, 1988.

\bibitem{BaCr} J.~Basto-Gon\c{c}alves and I.~Cruz, Analytic
  $k$-linearizability of some resonant Poisson structures,
  \emph{Lett.~Math.~Phys.~}\textbf{49}  (1999), 59--66.

\bibitem{Bochner} S.~Bochner, Compact groups of differentiable
  transformations, \emph{Ann.~of Math.~}\textbf{46} (1945), 372--381.

\bibitem{Bordemann} M.~Bordemann, A.~ Makhlouf and T.~Petit, 
D\'eformation par quantification et rigidit\'e des alg\`ebres
enveloppantes, preprint \emph{math.RA/0211416} (2002).

\bibitem{Bra} O.~Brahic, Normal forms of Poisson structures near a 
symplectic leaf, preprint \emph{math.SG/0403136} (2004).

\bibitem{CairGhy} G.~Cairns and E.~Ghys, The local linearization problem for 
smooth ${\rm SL}(n)$-actions, \emph{Enseign.~Math.~(2)}\textbf{43} (1997),
133--171. 

\bibitem{CaWe} A.~Cannas da Silva and A.~Weinstein, \emph{Geometric
Models for Noncommutative Algebras}, Berkeley Mathematics
Lectures, vol.~\textbf{10}, American Math.~Soc.~, Providence,
1999.

\bibitem{Cartan} H.~Cartan, Sur les groupes de transformation
  analytiques, \emph{Actualit\'es Scientifiques et Industrielles}
  \textbf{198} (1935).

\bibitem{Conn1} J.~Conn, Normal forms for analytic Poisson structures,
\emph{Annals of Math.}~\textbf{119} (1984), 576--601.

\bibitem{Conn2} J.~Conn, Normal forms for smooth Poisson structures,
\emph{Annals of Math.}~\textbf{121} (1985), 565--593.

\bibitem{Cr} M.~Crainic, Differentiable and algebroid cohomology, Van
Est isomorphisms, and characteristic classes, to appear in
\emph{Comment.~Math.~Helv.~}(preprint \emph{math.DG/0008064}).

\bibitem{CrFe1} M.~Crainic and R.~L.~Fernandes, Integrability of Lie
  brackets, \emph{Ann.~of Math.~(2)} \textbf{157}  (2003), 575--620.

\bibitem{CrFe2} M.~Crainic and R.~L.~Fernandes, Rigidity of Poisson
  brackets, in preparation.

\bibitem{CrFeMo} M.~Crainic, R.~L.~Fernandes and I.~Moerdijk,
  Deformation Cohomology, in preparation.

\bibitem{Desou} N.~Desolneux-Moulis, Lin\'earisation de certaines
  structures de Poisson de classe $C^\infty$. S\'eminaire de
  g\'eom\'etrie, 1985-1986, \emph{Publ.~D\'ep.~Math. Nouvelle
  S\'er. B}, \textbf{86-4}, Univ. Claude-Bernard, Lyon (1986), 55-68.

\bibitem{Dufour:Poisson} J.-P.~Dufour, Lin\'{e}arisation de certaines
  structures de Poisson, \emph{J.~Differential ~Geometry}~\textbf{32} (1990),
  no.~2, 415--428.

\bibitem{Dufour:algebroids} J.-P.~Dufour, Normal forms of Lie algebroids,
\emph{Banach Center Publications} \textbf{54} (2001), 35--41.

\bibitem{JP-Zung} J.~P.~Dufour and N.~T.~Zung, Nondegeneracy of the Lie
  algebra $\mathfrak{aff}(n)$ \emph{C.~R.~Acad.~Sci.~Paris
  S\'er.~I Math.~}\textbf{335} (2002), 1043--1046.

\bibitem{Book:JP-Zung} J.-P.~Dufour and N.~T.~Zung, \emph{Poisson Structures
and their Normal Forms}, book in preparation.

\bibitem{DuKo} J.~Duistermaat and J.~Kolk, \emph{Lie Groups},
  Springer-Verlag Berlin Heidelberg, 2000.

\bibitem{Fer} R.~L.~Fernandes, Lie algebroids, holonomy and
characteristic class, \emph{Adv.~in Math.~}\textbf{170} (2002),
119--179.

\bibitem{Flato} M.~Flato, G.~Pinczon and J.~Simon, Non linear 
representations of Lie groups, \emph{Ann. Sci. Ec. Norm. Sup}
\textbf{10} (1977), 405--418.

\bibitem{GuiSter} V.~Guillemnin and S.~Sternberg, Remarks on a paper
  of Hermann, \emph{Trans.~Amer.~Math. Soc.} \textbf{130} (1968),
  110--116.

\bibitem{Hamilton-NashMoser1982}
R.~Hamilton, The inverse function theorem of Nash and
Moser, \emph{Bull. Amer. Math. Soc.} (N.S.) \textbf{7} (1982),
no.~1, 65--222.

\bibitem{Hermann} R.~Hermann, The formal linearization of a semisimple
  Lie algebra of vector fields about a singular point,
  \emph{Trans.~Amer.~Math.~Soc.~}\textbf{130} (1968), 105--109.

\bibitem{Kushnirenko} A.~G.~Kushnirenko, Linear-equivalent action 
of a semi-simple Lie group in the neighborhood of a stationary point,
\emph{Funkts. Anal. Prilozh.} \textbf{1} (1967), 103--104. 

\bibitem{MoMr} I.~Moerdijk and J.~Mr\v{c}un, On integrability of
infinitesimal actions, \emph{Amer.~J.~Math.~}\textbf{124} (2002),
567--593.

\bibitem{Mol}  J.-C.~Molinier, \emph{Lin\'{e}arisation de structures de
Poisson}, Ph.D. Thesis, Montpellier 1993.

\bibitem{MoZung} Ph.~Monnier and N.T.~Zung, Levi decomposition for smooth
Poisson structures, preprint \emph{math.DG/0209004}.

\bibitem{Moser} J.~Moser, A new technique for the construction of
  solutions of nonlinear differential equations,
  \emph{Proc.~Nat.~Acad.~Sci.~USA} \textbf{47} (1961), 1824--1831.

\bibitem{Mr} J.~Mr\v{c}un, An extension of the Reeb stability theorem,
  \emph{Topology Appl.~}\textbf{70} (1996), 25--55. 

\bibitem{Wade} A.~Wade, Normalisation formelle de structures de Poisson,
\emph{C. R. Acad. Sci. Paris} S\'er. I Math. \textbf{324} (1997),
no.~5, 531--536.

\bibitem{Weibel} C.~Weibel, \emph{An Introduction to Homological Algebra},
  Cambridge University Press, Cambridge, 1994.

\bibitem{Wein} A.~Weinstein, The local structure of Poisson manifolds,
\emph{J.~Differential ~Geometry}~\textbf{18} (1983), 523--557.

\bibitem{Wein:rrank2} A.~Weinstein, Poisson geometry of the principal
series and non-linearizable structures, \emph{J.~Differential
~Geometry}~\textbf{25} (1987), 55--73.

\bibitem{Wein:linearization} A.~Weinstein, Linearization problem for
Lie algebroids and Lie groupoids, \emph{Lett. Math. Phys.}~\textbf{52}
(2000), no.~1, 93--102.

\bibitem{Wein:proper} A.~Weinstein, Linearization of regular proper
  groupoids, \emph{J.~Inst.~Math.~Jussieu} \textbf{1} (2002), 493--511.

\bibitem{Zung} N.T.~Zung, Levi decomposition of analytic Poisson
  structures and Lie algebroids, \emph{Topology}~\textbf{42} (2003),
  no.~6, 1403--1420.

\bibitem{Zung2} N.T.~Zung, A geometric proof of Conn's linearization
  theorem for analytic Poisson structures, preprint \emph{SG/0207263}.

\bibitem{Zung3} N.T.~Zung, Linearization of proper groupoids,
  preprint \emph{math.DG/0301297}.
\end{thebibliography}
\def\lllll{}

\end{document}